\input amstex
\documentstyle{amsppt}
\document
\topmatter
\title
The real section conjecture and Smith's fixed point theorem
\endtitle
\title
The real section conjecture and Smith's fixed point theorem for pro-spaces
\endtitle
\author
Ambrus P\'al
\endauthor
\date
September 12, 2009.
\enddate
\address
Department of Mathematics, 180 Queen's Gate, Imperial College,
London SW7 2AZ, United Kingdom\endaddress
\email a.pal\@imperial.ac.uk\endemail
\abstract We prove a topological version of the section conjecture for the profinite completion of the fundamental group of finite CW-complexes equipped with the action of a group of prime order $p$ whose $p$-torsion cohomology can be killed by finite covers. As an application we derive the section conjecture for the real points of a large class of varieties defined over the field of real numbers and the natural analogue of the section conjecture for fixed points of finite group actions on projective curves of positive genus defined over the field of complex numbers. 
\endabstract
\footnote" "{\it 2000 Mathematics Subject Classification. \rm 14P25, 14H30, 14G05.}
\endtopmatter

\heading 1. Introduction
\endheading

\definition{Notation 1.1} Let $X$ be a smooth, geometrically connected projective variety defined over $\Bbb R$. Let $\overline X$ and $\overline{\eta}$ be the base change of $X$ to $\Bbb C$ and a geometric point of $X$, respectively. Then Grothendieck's short exact sequence of \'etale fundamental groups for $X$ is:
$$\CD1@>>>\widehat\pi_1(\overline X,\overline{\eta})@>>>
\widehat\pi_1(X,\overline{\eta})@>>>\text{Gal}(\Bbb C|\Bbb R)@>>>1,\endCD\tag1.1.1$$
which is an exact sequence of profinite groups in the category of topological groups. Let $\Sigma$ be a set of primes and let $\Cal C_{\Sigma}$ denote the complete class of finite groups whose order is only divisible by primes in $\Sigma$. Let $\widehat\pi^{\Sigma}_1(X,\overline{\eta})$ and $\widehat\pi^{\Sigma}_1(\overline X,\overline{\eta})$ denote the maximal pro-$\Cal C_{\Sigma}$ quotients of the groups $\widehat\pi_1(X,\overline{\eta})$ and $\widehat\pi_1(\overline X,\overline{\eta})$, respectively. Assume that $2\in\Sigma$ and let
$$\CD1@>>>\widehat\pi_1^{\Sigma}(\overline X,\overline{\eta})
@>>>\widehat\pi_1^{\Sigma}(X,\overline{\eta})@>>>\text{Gal}(\Bbb C|\Bbb R)@>>>1\endCD\tag1.1.2$$
be the short exact sequence that we get from (1.1.1) by taking the maximal pro-$\Cal C_{\Sigma}$ quotients. Every $\Bbb R$-rational point $P\in X(\Bbb R)$ induces a section $\text{Gal}(\Bbb C|\Bbb R)\rightarrow\widehat\pi_1^{\Sigma}(X,\overline{\eta})$ of the sequence (1.1.2), well-defined up to conjugation. Let $\text{Sec}^{\Sigma}(X)$ denote the set of conjugacy classes of sections of (1.1.2) (in the category of profinite groups where morphisms are continuous homomorphisms). Then we have a map:
$$\iota^{\Sigma}_X:X(\Bbb R)\rightarrow\text{Sec}^{\Sigma}(X)$$
which sends every point $P\in X(\Bbb R)$ to the corresponding conjugacy class of sections. The set $X(\Bbb R)$ is a (possibly empty) differentiable manifold by the implicit function theorem. It is easy to prove that the map $\iota^{\Sigma}_X$ is constant on the connected components of $X(\Bbb R)$ (see for example Lemma 2.2 below). Let $\pi_0(X(\Bbb R))$ denote the set of connected components of $X(\Bbb R)$ and let
$$\iota^{\Sigma}_{X,0}:\pi_0(X(\Bbb R))\rightarrow\text{Sec}^{\Sigma}(X)$$
be the map which sends every connected component $C$ to the class $\iota^{\Sigma}_X(y)$ where $y\in C$ is arbitrary.
\enddefinition
\definition{Definition 1.2} We say that a connected, locally path-connected topological space $Y$ is $(\Sigma,p)$-erasable if for every $p$-torsion locally constant sheaf $\Cal F$ with finite stalks on $Y$ whose monodromy group is in $\Cal C_{\Sigma}$, for every positive integer $k$, and for every cohomology class $c\in H^k(Y,\Cal F)$ there is a connected finite topological cover $\pi:Z\rightarrow Y$ such that $\pi^*(c)\in H^k(Z,\pi^*(\Cal F))$ is zero, and for any point $z\in Z$ the image of the homomorphism $\pi_*:\pi_1(Z,z)\rightarrow\pi_1(Y,\pi(z))$ is a normal subgroup and the cokernel of $\pi_*$ is an element of $\Cal C_{\Sigma}$.
\enddefinition
 The motivation for this short note was to give a simple topological proof of the following:
\proclaim{Theorem 1.3} Assume that the differentiable manifold $X(\Bbb C)$ is $(\Sigma,2)$-erasable. Then the map $\iota^{\Sigma}_{X,0}$ is a bijection.
\endproclaim
The assumption in the theorem above is quite general: for example for a suitable choice of $\Sigma$ it is satisfied by varieties whose base change to $\Bbb C$ has a finite \'etale cover by a product of hyperbolic curves and abelian varieties, but there are other interesting examples, for example smooth non-isotrivial families of abelian varieties over hyperbolic curves. The result was only known in special cases previously (and were stated only for $\Sigma=\{\text{all primes}\}$). It was proved in the special case when $X$ is a hyperbolic curve or an abelian variety by Shinichi Mochizuki in [8], following an idea of Johan Huisman, using tools from the homotopy theory of algebraic varieties, namely a deep theorem of Cox (see [5]) and the theory developed in [9]. Recently another algebraic proof for hyperbolic curves was found by Jakob Stix, using a classical result of Witt on Brauer groups of real algebraic varieties and a bit of Kummer theory for the Jacobian (see [12]). Because the claim was formulated in analogy with Grothendieck's section conjecture in arithmetic geometry it is natural that the first proofs discovered were based on methods of arithmetic algebraic geometry itself.

On the other hand the result itself is essentially topological in nature. In fact the aim of this paper is to prove a general theorem in topology (Theorem 2.4) which might be considered a topological version of Theorem 1.3 and in particular it immediately implies the result above. This has the advantage that it trivially applies to smooth real varieties satisfying the same topological conditions, for example to non-algebraic real tori. On the other hand the argument is short, self-contained, and contains very little machinery compared to the existing proofs. In fact what we do is the following; first we show Theorem 4.5, an analogue of the classical Smith fixed point theorem for pro-spaces, using standard tools from algebraic topology. Then the latter is applied to a certain pro-object in the category of pointed topological spaces, which serves as a universal covering space for (certain) finite covers, to derive Theorem 2.4. Hence the argument also realizes the original philosophy of Grothendieck (see [6]): the section conjecture is a fixed point theorem for the pro-universal cover. Of course this strategy is also inspired by the work of Artin-Mazur. Our main result is general enough to enable us to derive also with ease an analogue of the section conjecture for arbitrary finite group actions on smooth projective algebraic curves of positive genus defined over $\Bbb C$ (Theorem 5.11).
\definition{Contents 1.4} In the next chapter we formulate the main result (Theorem 2.4) of this paper. In the third chapter for every complete class $\Cal C$ we construct a pro-object in the category of pointed topological spaces that will serve as a universal covering space for finite Galois covers of a given pointed connected locally contractible topological space $X$, and under certain assumptions on $X$ and $\Cal C$ we show that it is cohomologically acyclic mod $p$. We prove a generalization of the Smith fixed point theorem for pro-objects in the fourth chapter following the argument due to A. Borel, and with its help and the results of the previous chapter we conclude the proof of Theorem 2.4. In the last chapter first we present a mild generalization of the theory of good groups, and use these results to provide examples satisfying the conditions of Definition 1.2. Then we use our main result to derive Theorem 1.3 and a similar application to arbitrary finite group actions on smooth projective algebraic curves of positive genus. 
\enddefinition
\definition{Acknowledgement 1.5} I wish to thank Akio Tamagawa for some useful discussion concerning the contents of the last chapter. The author was partially supported by the EPSRC grant P19164. 
\enddefinition

\heading 2. Basic setup
\endheading

\definition{Definition 2.1} In this paper by a topological space we always mean a Hausdorff topological space. For every topological space let Aut$(Y)$ denote its group of self-homeomorphisms. Let $X$ be a connected locally contractible topological space. Let $G$ be a finite group and assume that a non-trivial continuous left action $h:G\times X\rightarrow X$ of $G$ on $X$ is given. Let $x$ be a point of $X$ and let $(\widetilde X,\widetilde x)$ be the universal cover of the pointed space $(X,x)$. Let $\pi:\widetilde X\rightarrow X$ be the covering map (with the property $\pi(\widetilde x)=x$). We define $\pi_1(X,G,x)$ to be the following subgroup of Aut$(\widetilde X)\times G$:
$$\pi_1(X,G,x)=\{(\phi,g)\in\text{Aut}(\widetilde X)\times G|
\pi\circ\phi=h(g)\circ\pi\}.$$
By definition the group $\pi_1(X,G,x)$ acts on $\widetilde X$. Moreover the fundamental group $\pi_1(X,x)$, considered as the group of deck transformations of the cover $\pi:\widetilde X\rightarrow X$, is a subgroup of $\pi_1(X,G,x)$ via the map $\phi\mapsto(\phi,1)$. Because every continuous map $X\rightarrow X$ can be lifted by the universal property of the covering map $\pi$ we have an exact sequence:
$$\CD1@>>>\pi_1(X,x)@>>>\pi_1(X,G,x)@>\psi_{X,G}>>G@>>>1.\endCD\tag2.1.1$$
Let $X^G$ denote the set of fixed points of the action of $G$ on $X$. Let $y\in X^G$ and choose a $\widetilde y\in\widetilde X$ such that $\pi(\widetilde y)=y$. For every $g\in G$ there is a unique $h_{\widetilde y}(g)\in\pi_1(X,G,x)$ which fixes $\widetilde y$ and $\psi_{X,G}(h_{\widetilde y}(g))=g$. The map $h_{\widetilde y}:G\rightarrow\pi_1(X,G,x)$ is a homomorphism which is also a section of the exact sequence (2.1.1) by construction. Let $\widetilde y'\in\widetilde X$ be another point such that $\pi(\widetilde y)=y$. Then there is a unique $\gamma\in\pi_1(X,x)$ such that $\widetilde y'=\gamma(\widetilde y)$. Clearly $h_{\widetilde y'}=\gamma^{-1}h_{\widetilde y}\gamma$.

We say that two sections $s_1:G\rightarrow\pi_1(X,G,x)$ and $s_2:G\rightarrow\pi_1(X,G,x)$ of the exact sequence (2.1.1) are conjugate if there a $\gamma\in\pi_1(X,x)$ such that $s_1=\gamma^{-1}s_2\gamma$. Obviously the latter is an equivalence relation. As usual we will call the equivalence classes of this equivalence relation conjugacy classes. Let $\text{Sec}(X,G)$ denote the set of conjugacy classes of sections of (2.1.1). By the above the conjugacy class of $h_{\widetilde y}$ does not depend on the choice of the lift $\widetilde y$. Hence we get a well-defined map:
$$\iota_X:X^G\rightarrow\text{Sec}(X,G)$$
which associates to $y$ the conjugacy class of $h_{\widetilde y}$. 
\enddefinition
\proclaim{Lemma 2.2} The map $\iota_X$ is constant on the path-connected components of $X^G$.
\endproclaim
\definition{Proof} Let $C\subseteq X^G$ be a path-connected component and let $y$, $z\in C$. Let $p:[0,1]\rightarrow C$ be a continuous path with the property $p(0)=y$ and $p(1)=z$. Choose a $\widetilde y\in\widetilde X$ such that $\pi(\widetilde y)=y$. As $\widetilde X\rightarrow X$ is a topological covering there is a unique lift $\widetilde p:[0,1]\rightarrow\widetilde X$ of the path $p$ with base point $\widetilde y$. (The latter means that $\widetilde p$ is a continuous map such that $\pi\circ\widetilde p=p$ and $\widetilde p(0)=\widetilde y$.) For every $g\in G$ let $\widetilde p^g$ be the image of $\widetilde p$ under $h_{\widetilde y}(g)$. Then the set $I_g=\{t\in[0,1]|\widetilde p(t)=\widetilde p^g(t)\}$ is closed. But the covering map $\widetilde X\rightarrow X$ is locally a homeomorphism and $p$ is fixed by $G$ hence this set is also open. As $0\in I_g$ the set $I_g$ is not empty. Because $[0,1]$ is connected as a topological space we get that $I_g=[0,1]$ for every $g\in G$. Therefore the path $\widetilde p$ lies in the set $\widetilde X^G$ of fixed points of the action of $G$ on $\widetilde X$ given by $h_{\widetilde y}$. In particular the endpoint $\widetilde z=\widetilde p(1)$ of $\widetilde p$ is also an element of $\widetilde X^G$. By definition $\pi(\widetilde z)=z$. Because $h_{\widetilde z}(g)\in\pi_1(X,G,x)$ is the unique lift of $h(g)$ which fixes $\widetilde z$ we get that $h_{\widetilde z}(g)=h_{\widetilde y}(g)$ for every $g\in G$.\ $\square$
\enddefinition
\definition{Definition 2.3} Let $\Sigma$ be a set of primes and let $\Cal C_{\Sigma}$ denote the complete class introduced in Notation 1.1. For every group $H$ let $H^{\Sigma}$ denote the $\Cal C_{\Sigma}$-completion of $H$ (for the definition of the latter see page 26 of [2]). Moreover let $\pi^{\Sigma}_1(X,x)$ and $\pi^{\Sigma}_1(X,G,x)$ denote the profinite groups $\pi_1(X,x)^{\Sigma}$ and $\pi_1(X,G,x)^{\Sigma}$, respectively. By applying the $\Cal C_{\Sigma}$-completion functor to (2.1.1) we get the following exact sequence:
$$\CD1@>>>\pi^{\Sigma}_1(X,x)@>>>
\pi^{\Sigma}_1(X,G,x)@>>>G@>>>1\endCD\tag2.3.1$$ 
when $G$ is in $\Cal C_{\Sigma}$. We say that two sections $s_1:G\rightarrow\pi^{\Sigma}_1(X,G,x)$ and $s_2:G\rightarrow\pi^{\Sigma}_1(X,G,x)$ of the exact sequence (2.3.1) are conjugate if there a $\gamma\in\pi^{\Sigma}_1(X,x)$ such that $s_1=\gamma^{-1}s_2\gamma$. This is again an equivalence relation and we will call the equivalence classes of this equivalence relation conjugacy classes, too. Let $\text{Sec}^{\Sigma}(X,G)$ denote the set of conjugacy classes of sections of (2.3.1). The $\Cal C_{\Sigma}$-completion induces a natural map $f^{\Sigma}:\text{Sec}(X,G)\rightarrow\text{Sec}^{\Sigma}(X,G)$. For every topological space $M$ let $\pi_0(M)$ denote the set of connected components of $M$. By the above when $X^G$ is locally path-connected there is a well-defined map:
$$\iota^{\Sigma}_{X,0}:\pi_0(X^G)\rightarrow
\text{Sec}^{\Sigma}(X,G)\tag2.3.2$$
which sends every connected component $C$ to the class $f^{\Sigma}(\iota_X(y))$ where $y\in C$ is arbitrary.
\enddefinition
Assume now that $X$ is a $G$-CW complex. Then $X^G$ is also a CW-complex which we get by attaching only those $G$-cells of $X$ which are of the type $G/G$ (see Definition 1.1.1 of [1] on page 3) hence the map $\iota^{\Sigma}_{X,0}$ of $(2.3.2)$ is well-defined. Our main result is the following:
\proclaim{Theorem 2.4} Assume that $X$ is a finite $G$-CW complex and $G$ is a finite group of prime order $p$. Also suppose that $p\in\Sigma$ and $X$ is $(\Sigma,p)$-erasable. Then the map $\iota^{\Sigma}_{X,0}$ is a bijection.
\endproclaim
This result will be proved in the next two sections. 

\heading 3. The pro-universal cover
\endheading

\definition{Notation 3.1} We will follow the terminology of the book [2]. For every category $C$ let pro-$C$ denote the category of pro-objects of $C$ and let $C^{\text{op}}$ denote the opposite of $C$, respectively. Let $\Cal T$ denote the category whose objects are topological spaces and whose morphisms are continuous maps between them. Let $\Cal T_0$ denote the category whose objects are connected pointed locally contractible topological spaces and whose morphisms are pointed continuous maps between them. There is a faithful forgetful functor $\Cal T_0\rightarrow\Cal T$. We will denote the image of every object or pro-object $X$ of $\Cal T_0$ by the same symbol $X$ under this functor by slight abuse of notation. For every $k\in\Bbb N$, abelian group $A$ and pro-object $X=\{X_i\}_{i\in I}$ of the category $\Cal T$ let
$$H^k(X,A)=\varinjlim_{i\in I}H^k(X_i,A)$$
be the inductive limit of the groups $H^k(X_i,A)$. The assignment $X\mapsto H^k(X,A)$ is clearly a contravariant functor from pro-$\Cal T$ to the category of abelian groups. As usual for every morphism $\pi:X\rightarrow Y$ the symbol $\pi^*$ will denote the homomorphism $H^k(Y,A)\rightarrow H^k(X,A)$ induced by functoriality. Fix a prime number $p$ and for every object $X$ of pro-$\Cal T$ let $H^k(X)$ denote the group $H^k(X,\Bbb F_p)$. 
\enddefinition
\definition{Definition 3.2} Let $X$ be an object of $\Cal T_0$ and let $\Cal C$ be a complete class of finite groups. Let $I(X,\Cal C)$ be the set of all normal subgroups $K$ of $\pi_1(X)$ of finite index such that the quotient group $\pi_1(X)/K$ is in $\Cal C$. Let $\Cal I(X,\Cal C)$ be the category whose objects are the elements of $I(X,\Cal C)$ and for every pair of objects $K$, $L\in I(X,\Cal C)$ the set of morphisms from $K$ to $L$ consists of the ordered pair $\phi_{K,L}=(K,L)$, if $L\leq K$, and it is empty, otherwise. Clearly $\Cal I(X,\Cal C)$ is a small filtering index category. For every $K\in I(X,\Cal C)$ let $\pi_K:X_K\rightarrow X$ be the covering map of pointed topological spaces corresponding to finite index subgroup $K\leq\pi_1(X)$. For every pair of $K$, $L\in I(X,\Cal C)$ such that $L\leq K$ there is a unique map $\pi_{K,L}:X_L\rightarrow X_K$ such that $\pi_L=\pi_K\circ\pi_{K,L}$. Therefore there is a functor $\Cal F_{\Cal C}:\Cal I(X,\Cal C)^{\text{op}}\rightarrow\Cal T_0$ which assigns to every object $K$ the pointed space $X_K$ and to every morphism $\phi_{K,L}$ the map $\pi_{K,L}$. When $\Cal C=\Cal C_{\Sigma}$ for some set of primes $\Sigma$ let $\widetilde X_{\Sigma}$ denote the object of pro-$\Cal T_0$ corresponding to the functor $\Cal F_{\Cal C_{\Sigma}}$. In this case  the set of open normal subgroups of $\pi_1(X)^{\Sigma}$ and $I(X,\Cal C)$ are in a natural one to one and onto correspondence, and we will not distinguish them in all that follows.
\enddefinition 
We say that an object $X$ of pro-$\Cal T$ is cohomologically acyclic mod $p$ if 
$$H^k(X)=\cases\Bbb F_p,&\text{if $k=0$,}\\
0,&\text{otherwise.}\endcases$$
\proclaim{Proposition 3.3} Let $\Sigma$ be a set of primes such that $p\in\Sigma$ and let $X$ be a pointed connected finite CW-complex which is $(\Sigma,p)$-erasable. Then $\widetilde X_{\Sigma}$ is cohomologically acyclic mod $p$.
\endproclaim
\definition{Proof} First note that for every $K\in I(X,\Cal C_{\Sigma})$ the space $X_K$ is connected, hence $H^0(X_K)=\Bbb F_p$. By taking the limit we get that $H^0(\widetilde X_{\Sigma})=\Bbb F_p$, too. Now it will be sufficient to prove that for every positive integer $k$, for every $K\in I(X,\Cal C_{\Sigma})$ and for every $c\in H^k(X_K)$ there is an $L\in I(X,\Cal C_{\Sigma})$ such that $L\leq K$ and $\pi_{K,L}^*(c)=0$. Hence the claim follows from Lemma 3.4 below.\ $\square$
\enddefinition
\proclaim{Lemma 3.4} For every $K\in I(X,\Cal C_{\Sigma})$ the space $X_K$ is $(\Sigma,p)$-erasable.
\endproclaim
\definition{Proof} Let $\Cal F$ be a $p$-torsion locally constant sheaf $\Cal F$ with finite stalks on $X_K$; then the direct image sheaf $(\pi_K)_*(\Cal F)$ is also a $p$-torsion locally constant sheaf with finite stalks on $X$. Assume moreover that the monodromy group of $\Cal F$ is in $\Cal C_{\Sigma}$. Then the monodromy group of $(\pi_K)_*(\Cal F)$ is also in $\Cal C_K$, because it is the extension of a subgroup of $\pi_1(X)/K$ by the monodromy group of $\Cal F$ and $\Cal C_{\Sigma}$ is a complete class.  Now let $k$ be a positive integer and let $c\in H^k(X_K,\Cal F)$ be a cohomology class. For every $L\in I(X,\Cal C_{\Sigma})$ such that $L\leq K$ let $\text{Res}_L^1:H^k(X,(\pi_K)_*(\Cal F))\rightarrow H^k(X_L,(\pi_L)^*(\pi_K)_*(\Cal F))$ and $\text{Res}^2_L:H^k(X_K,\Cal F)\rightarrow H^k(X_L,(\pi_{K,L})^*(\Cal F))$ be the natural restriction maps. Moreover for every $L$ as above let
$$\text{Pr}_L:H^k(X_L,(\pi_L)^*(\pi_K)_*(\Cal F))\rightarrow H^k(X_L,(\pi_{K,L})^*(\Cal F))$$
be the map induced by the natural projection:
$$(\pi_L)^*(\pi_K)_*(\Cal F)=
(\pi_{K,L})^*(\pi_K)^*(\pi_K)_*(\Cal F)
\rightarrow(\pi_{K,L})^*(\Cal F).$$
Because the composition $\text{Pr}_K\circ\text{Res}_K^1:H^k(X,(\pi_K)_*(\Cal F))\rightarrow H^k(X_K,\Cal F)$ is an isomorphism there is a $d\in H^k(X,(\pi_K)_*(\Cal F))$ such that $\text{Pr}_K(\text{Res}_K^1(d))=c$. By assumption there is an $L\in I(X,\Cal C_{\Sigma})$ such that $L\leq K$ and $\text{Res}^1_L(d)=0$. Since
$$\text{Res}^2_L(c)=\text{Res}^2_L(\text{Pr}_K(\text{Res}_K^1(d)))=\text{Pr}_L((\text{Res}^1_L(d))=0,$$
the claim is now clear.\ $\square$
\enddefinition
\definition{Definition 3.5} Let $pt$ denote the topological space whose underlying set consists of one point. For every pro-object $X$ of $\Cal T$ we define its set of points $\underline X$ to be the set of morphisms from $pt$ to $X$ in the category pro-$\Cal T$. The assignment $X\mapsto\underline X$ is clearly a functor from pro-$\Cal T$ to the category of sets. For every morphism $\phi:X\rightarrow Z$ in the category pro-$\Cal T$ and $z\in\underline Z$ we say that $z$ lies in the image of $\phi$ if $z$ lies in the image of the corresponding map $\underline X\rightarrow\underline Z$. When $X$ is an object of $\Cal T$ we will continue to denote $\underline X$ by $X$, following the usual abuse of notation. 
\enddefinition
Let $\Cal I$ be a small filtering index category and let $I$ denote the set of objects of $\Cal I$. Let $\bold X:\Cal I^{\text{op}}\rightarrow\Cal T$ be a functor and let $X=\{\bold X(i)\}_{i\in I}$ be the corresponding pro-object of the category $\Cal T$.  The following claim is very well known. 
\proclaim{Lemma 3.6} Assume that for every $i\in I$ the topological space $\bold X(i)$ is compact. Then $\underline X$ is non-empty if and only if for every $i\in I$ the topological space $\bold X(i)$ is non-empty.
\endproclaim
\definition{Proof} The standard argument (see Proposition 8 of [4] on page 89, where the special case of directed sets is presented)  works here without much modification.\ $\square$
\enddefinition
\definition{Definition 3.7} Let $\Cal I,I,\bold X$ and $X$ be as above. Let $G$ be a finite group and assume that for every $i\in I$ a continuous left action $h_i:G\times\bold X(i)\rightarrow\bold X(i)$ of $G$ on $\bold X(i)$ is given such that for every morphism $\phi:i\rightarrow j$ of $\Cal I^{\text{op}}$ the map $\bold X(\phi):\bold X(i)\rightarrow\bold X(j)$ is equivariant with respect to $h_i$ and $h_j$. In this case we will say that the collection $h=\{h_i\}_{i\in I}$ is a continuous left action of $G$ on $X$. For every $i\in I$ let $\bold X^G(i)$ denote the subspace of fixed points of the action $h_i$ on $\bold X(i)$. Let $\bold X^G:\Cal I^{\text{op}}\rightarrow\Cal T$ be a functor which assigns to every $i\in I$ the topological space $\bold X^G(i)$ and to every morphism $\phi:i\rightarrow j$ of $\Cal I^{\text{op}}$ the restriction of $\bold X(\phi)$ onto $\bold X^G(i)$; this functor is indeed well-defined (that is, $\bold X(\phi)$ maps $\bold X^G(i)$ into $\bold X^G(j)$) by our assumptions. Let $X^G=\{\bold X^G(i)\}_{i\in I}$ be the corresponding pro-object of the category $\Cal T$. We will call the latter the subspace of fixed points of the action $h$. 
\enddefinition
\definition{Definition 3.8} Again let $X$ be a connected locally contractible topological space and let $G$ be a finite group and assume that a non-trivial continuous left action $h:G\times X\rightarrow X$ of $G$ on $X$ is given, as in Definition 2.1. Fix a base point $x$ of $X$ and by slight abuse of notation let $X$ also denote the pointed space $(X,x)$ when this does not cause confusion. Let $\Sigma$ be a set of primes and assume that $G\in\Cal C_{\Sigma}$. Let $\widetilde X$ denote the universal cover of the pointed space $X$ and for every $K\in I(X,\Cal C_{\Sigma})$ let $\pi^K:\widetilde X\rightarrow X_K$ be the covering map. Let $j:\pi_1(X,G,x)\rightarrow\pi_1^{\Sigma}(X,G,x)$ be the map furnished by the $\Cal C_{\Sigma}$-completion functor. Let $I(X,\Cal C_{\Sigma},G)$ denote the set of 
open subgroups $K\leq\pi^{\Sigma}_1(X,x)$ which are also normal as a subgroup of the larger group $\pi^{\Sigma}_1(X,G,x)$. Let $s:G\rightarrow\pi^{\Sigma}_1(X,G,x)$ be a section of (2.3.1), and 
for every $K\in I(X,\Cal C_{\Sigma},G)$ and for every $g\in G$ pick an element $\widetilde g_K\in\pi_1(X,G,x)$ such that $j(\widetilde g_K)$ and $s(g)$ are in the same right $K$-coset. By our assumptions on $K$ the element $\widetilde g_K$ normalizes the subgroup $j^{-1}(K)$ hence there is a unique homeomorphism $h^K_s(g):X_K\rightarrow X_K$ such that $\pi^K\circ\widetilde g_K=h^K_s(g)\circ\pi^K$. Moreover $h^K_s(g)$ is independent of the choice of $\widetilde g_K$ and the map $g\mapsto h^K_s(g)$ defines a continuous left-action of $G$ on $X_K$ which will be denoted simply by $h^K_s$. Let $\Cal I(X,\Cal C_{\Sigma},G)$ denote the full subcategory of $\Cal I(X,\Cal C_{\Sigma})$ whose set of objects is $I(X,\Cal C_{\Sigma},G)$. Then 
$\Cal I(X,\Cal C_{\Sigma},G)$ is a cofinal subcategory of $\Cal I(X,\Cal C_{\Sigma})$ therefore the collection $h_s=\{h^K_s\}_{K\in I(X,\Cal C_{\Sigma},G)}$ defines a continuous left action of $G$ on $\widetilde X_{\Sigma}$. Let $\widetilde X^G_{\Sigma}$ denote the subspace of fixed points of the action $h_s$. Note that there is a natural map $\widetilde X_{\Sigma}^G\rightarrow X^G$. 
\enddefinition
In the next two lemmas we will assume that $X^G$ is locally path-connected. 
\proclaim{Lemma 3.9} The image of the map $\widetilde X_{\Sigma}^G\rightarrow X^G$ consists of the union of some connected components of $X^G$.
\endproclaim
\definition{Proof} For the sake of simple notation for every $K\in I(X,\Cal C_{\Sigma},G)$ let $X^G_K$ denote the set of fixed points of the action $h^K_s$. Let $C\subseteq X^G$ be a connected component and let $y$, $z\in C$. We have to show that if $y$ lies in the image of the map $\widetilde X_{\Sigma}^G\rightarrow X^G$ then so does $z$. Hence assume that for every $K\in I(X,\Cal C_{\Sigma},G)$ there is a $y_K\in X_K^G$ such that $\pi_K(y_K)=y$ and these points are compatible in the sense that for every $K$, $L\in I(X,\Cal C_{\Sigma},G)$ such that $L\leq K$ we have $\pi_{K,L}(y_L)=y_K$. Let $p:[0,1]\rightarrow C$ be a continuous path with the property $p(0)=y$ and $p(1)=z$. For every $K\in I(X,\Cal C_{\Sigma},G)$ the map $\pi_K:X_K\rightarrow X$ is a topological covering hence there is a unique lift $p_K:[0,1]\rightarrow X_K$ of the path $p$ with base point $y_K$. These maps are compatible, that is for every $K$, $L\in I(X,\Cal C_{\Sigma},G)$ such that $L\leq K$ we have $\pi_{K,L}\circ p_L=p_K$. By repeating the argument in the proof of Lemma 2.2 we get that the path $p_K$ lies in $X^G_K$. In particular its endpoint $z_K=p_K(1)$ is also an element of $X_K^G$. These points are compatible hence the system $\{z_k\}_{K\in I(X,\Cal C_{\Sigma},G)}$ defines a point of $\widetilde X^G_{\Sigma}$ whose image under the map is $\widetilde X_{\Sigma}^G\rightarrow X^G$ is $z$.\ $\square$
\enddefinition
\proclaim{Lemma 3.10} The connected component $C\in\pi_0(X^G)$ is in the image of  the map $\widetilde X_{\Sigma}^G\rightarrow X^G$ if and only if the conjugacy class of $s$ is equal to $\iota^{\Sigma}_{X,0}(C)$.
\endproclaim
\definition{Proof} Let $y\in C$ and choose a $\widetilde y\in\widetilde X$ such that $\pi(\widetilde y)=y$ (where $\pi:\widetilde X\rightarrow X$ is the covering map). For every $g\in G$ let $h_{\widetilde y}(g)\in\pi_1(X,G,x)$ denote the unique map which fixes $\widetilde y$ and $\pi\circ h_{\widetilde y}(g)=h(g)\circ\pi$ as in Definition 2.1. For every $K\in I(X,\Cal C_{\Sigma},G)$ let $\widetilde y_K=\pi^K(\widetilde y)$ and for every such $K$ and $g\in G$ let $h^K_{\widetilde y}(g):X_K\rightarrow X_K$ be the unique homeomorphism such that $\pi^K\circ h_{\widetilde y}(g)=h^K_{\widetilde y}(g)\circ\pi^K$. First assume that $y$ lies in the image of the map $\widetilde X_{\Sigma}^G\rightarrow X^G$, that is, for every $K\in I(X,\Cal C_{\Sigma},G)$ there is a point $y_K\in X_K^G$ such that $\pi_K(y_K)=y$ and these points are compatible in the sense of the proof of Lemma 3.9 above. Then $\pi_K(\widetilde y_K)=\pi_K(y_K)$ hence there is a unique element $\gamma_K$ in the group $\pi_1^{\Sigma}(X)/K$ of deck transformations of the cover $\pi_K:X_K\rightarrow X$ such that $\gamma_K(y_K)=\widetilde y_K$. For every $K$, $L\in I(X,\Cal C_{\Sigma},G)$ such that $L\leq K$ the image of $\gamma_L$ with respect to the quotient map $\pi_1^{\Sigma}(X)/L\rightarrow \pi_1^{\Sigma}(X)/K$ is $\gamma_K$. Hence there is a unique $\gamma\in\pi_1^{\Sigma}(X)$ such that the image of $\gamma$ with respect to the quotient map $\pi_1^{\Sigma}(X)\rightarrow \pi_1^{\Sigma}(X)/K$ is $\gamma_K$ for every $K\in I(X,\Cal C_{\Sigma},G)$. Clearly 
$\gamma_K^{-1}h^K_{\widetilde y}(g)\gamma_K=h^K_s(g)$ for every $K\in I(X,\Cal C_{\Sigma},G)$ and $g\in G$ hence $\gamma^{-1}h_{\widetilde y}(g)\gamma=s(g)$ for every $g\in G$, too. Assume now the converse, that is, there is a $\gamma\in\pi_1^{\Sigma}(X)$ such that $\gamma^{-1}h_{\widetilde y}(g)\gamma=s(g)$ for every $g\in G$. Let $\gamma_K$ denote the image of $\gamma$ with respect to the quotient map $\pi_1^{\Sigma}(X)\rightarrow \pi_1^{\Sigma}(X)/K$ for every $K\in I(X,\Cal C_{\Sigma},G)$. Let $y_K=\gamma_K^{-1}(\widetilde y_K)$ for every $K$ as above; these points are compatible in the sense of the proof of Lemma 3.9 above, $\pi_K(y_K)=y$ and
$y_K\in X^G_K$. Hence $y$ lies in the image of  the map $\widetilde X_{\Sigma}^G\rightarrow X^G$.\ $\square$
\enddefinition

\heading 4. The Smith fixed point theorem
\endheading

\definition{Definition 4.1} Let $\Cal T_{pair}$ be the category whose objects are ordered pairs $(X,Y)$ where $X$ is a topological space and $Y$ is a closed subspace and a morphism from an object $(X,Y)$ to another object $(X',Y')$ is a continuous map $f:X\rightarrow X'$ such that $f(Y)\subseteq Y'$. A pro-object of $\Cal T_{pair}$ can be also described as follows; let $\Cal I,I,\bold X$ and $X$ be as in Definition 3.7. Assume that for every $i\in I$ a closed subspace $\bold Y(i)$ of $\bold X(i)$ is given such that for every morphism $\phi:i\rightarrow j$ of $\Cal I^{\text{op}}$ the image of $\bold Y(i)$ with respect to the map $\bold X(\phi)$ lies in $\bold Y(j)$. Then let $(\bold X,\bold Y):\Cal I^{\text{op}}\rightarrow\Cal T_{pair}$ be a functor which assigns to every $i\in I$ the pair $(\bold X(i),\bold Y(i))$ and to every morphism $\phi:i\rightarrow j$ of $\Cal I^{\text{op}}$ the map $\bold X(\phi)$; this functor is clearly well-defined. We will denote the pro-object of $\Cal T_{pair}$ corresponding to the functor $(\bold X,\bold Y)$ by $(X,Y)$ by slight abuse of notation, where $Y$ is the pro-object $\{\bold Y(i)\}_{i\in I}$ of $\Cal T$. 
\enddefinition
\definition{Definition 4.2} For every $k\in\Bbb N$, abelian group $A$ and object $(X,Y)$ in pro-$\Cal T_{pair}$ let
$$H^k(X,Y,A)=\varinjlim_{i\in I}H^k(X_i,Y_i,A)$$
be the inductive limit of the groups $H^k(X_i,Y_i,A)$ where $X=\{X_i\}_{i\in I}$ and $Y=\{Y_i\}_{i\in I}$. The assignment $(X,Y)\mapsto H^k(X,Y,A)$ is  a contravariant functor from pro-$\Cal T_{pair}$ to the category of abelian groups. For every morphism $\pi:(X,Y)\rightarrow(X',Y')$ in pro-$\Cal T_{pair}$ let the symbol $\pi^*$ denote the homomorphism $H^k(X',Y',A)\rightarrow H^k(X,Y,A)$ induced by functoriality. Because taking the direct limit is an exact functor there is a long exact sequence:
$$\CD\rightarrow H^{k-1}(Y,A)@>>>H^k(X,Y,A)@>>>
H^k(X,A)@>>>H^k(Y,A)\rightarrow\endCD$$
which we will call the cohomological long exact sequence of the pair $(X,Y)$.
The latter is also functorial. For every $(X,Y)$ as above let $H^k(X,Y)$ denote the group $H^k(X,Y,\Bbb F_p)$ where $p$ is the prime number we fixed in Definition 3.1. 
\enddefinition
\definition{Definition 4.3} Let $\Cal I,I,\bold X$ and $X$ be as in Definition
4.1, let $G$ be a finite group, and let $h=\{h_i\}_{i\in I}$ be a continuous left action of $G$ on $X$. For every $i\in I$ let $\bold X/G(i)$ denote the quotient of $\bold X(i)$ by $G$ with respect to the action $h_i$ and let
$\pi_i:\bold X(i)\rightarrow\bold X/G(i)$ be the quotient map. There is a unique functor $\bold X/G:\Cal I^{\text{op}}\rightarrow\Cal T$ which assigns to every $i\in I$ the topological space $\bold X/G(i)$ and to every morphism $\phi:i\rightarrow j$ of $\Cal I^{\text{op}}$ the unique map $\bold X/G(\phi):\bold X/G(i)\rightarrow\bold X/G(j)$ such that
$\pi_j\circ\bold X(\phi)=\bold X/G(\phi)\circ\pi_i$. Let $X/G=\{\bold X/G(i)\}_{i\in I}$ be the corresponding pro-object of the category $\Cal T$; we will call it the quotient of $X$ by the action $h$. Then the collection $\pi=\{\pi_i\}_{i\in I}$ defines a morphism from $X$ to $X/G$ which we will call the quotient map. 
\enddefinition
Let $(X,Y)$ be an object in pro-$\Cal T_{pair}$ and assume that $h=\{h_i\}_{i\in I}$ is a continuous left action of the finite group $G$ on $X$ such that for every $i\in I$ the action $h_i$ leaves $Y_i\subseteq X_i$ invariant, where $X=\{X_i\}_{i\in I}$ and $Y=\{Y_i\}_{i\in I}$. Then the action $h_i$ induces an action of $G$ on the cohomology group $H^i(X_i,Y_i,A)$ for every $i\in I$. Because the system $\{h_i\}_{i\in I}$ is a continuous left action of $G$ on $X$ the transition maps between the $G$-modules $H^i(X_i,Y_i,A)$ are $G$-module homomorphisms, so the inductive limit $H^k(X,Y,A)$ is also
equipped with the structure of a $G$-module. Assume now that for every $i\in I$ the topological space $X_i$ is compact and the action $h_i$ is free on $X_i-Y_i$.
\proclaim{Lemma 4.4} There is a spectral  sequence $H^j(G,H^k(X,Y))\Rightarrow H^{j+k}(X/G,Y/G)$.
\endproclaim
We will call the spectral  sequence above the Hochschild-Serre spectral sequence of $(X,Y)$. 
\definition{Proof} By 1.1 of [3] on page 29 for every $i\in I$ there is a spectral sequence $S_i$:
$$H^j(G,H^k(X_i,Y_i,A))\Rightarrow H^{j+k}(X_i/G,Y_i/G,A).$$
Moreover the assignment $i\mapsto S_i$ is a functor from the indexing category of the pro-object $(X,Y)$ into the category of spectral sequences. Because taking the direct limit is an exact functor the claim is now clear.\ $\square$
\enddefinition
\proclaim{Theorem 4.5} Let $G$ be a finite group of order $p$, let $X=\{X_i\}_{i\in I}$ be a pro-object of the category $\Cal T$ which is cohomologically acyclic modulo $p$ and let $h=\{h_i\}_{i\in I}$ be a continuous left action of $G$ on $X$. Assume that $X_i$ is a finite $G$-CW complex with respect to the action $h_i$ for every $i\in I$. Also assume that there is a natural number $N$ such that $X_i$ has no cells of dimension bigger than $N$ for every $i\in I$. Then the subspace $X^G$ of fixed points of $h$ is cohomologically acyclic mod $p$, too. 
\endproclaim
\definition{Proof} Let $Y$ denote the quotient of $X$ by the action $h$ and let $\pi=\{\pi_i\}_{i\in I}$ be the quotient map $\pi:X\rightarrow Y$. The restriction of $\pi_i$ onto $X_i^G$ is a topological embedding for every $i\in I$.  We'll denote the pro-object $\{\pi_i(X^G_i)\}_{i\in I}$ also by the symbol $X^G$ by slight abuse of notation. Because $G$ has prime order,
it acts freely on $X_i-X^G_i$ for every $i\in I$. Hence there is a Hochschild-Serre spectral sequence:
$$H^j(G,H^k(X,X^G))\Rightarrow H^{j+k}(Y,X^G)$$
by Lemma 4.4.
\enddefinition
\proclaim{Lemma 4.6} The following holds:
\roster
\item"$(i)$" we have:
$$H^k(X,X^G)=H^k(Y,X^G)=0$$
for every $k>N+1$,
\item"$(ii)$" the group $H^k(X,X^G)$ is equal to the image of $\pi^*$ for every $k\geq0$,
\item"$(iii)$" the term $E_2^{j,k}$ of the Hochschild-Serre spectral sequence of $(X,X^G)$ is isomorphic to $H^k(X,X^G)$ and is formed of permanent cocycles for every $k\geq0$.
\endroster
\endproclaim
\definition{Proof} In order to prove claim $(i)$ it will be sufficient to show that $H^k(X)=H^k(Y)=H^k(X^G)=0$ for every $k>N$ because of the existence of cohomological long exact sequences for the pairs $(X,X^G)$ and $(Y,X^G)$. In order to do so, it will be sufficient to show that $H^k(X_i)=H^k(Y_i)=H^k(X^G_i)=0$ for every $k>N$ and $i\in I$. But the latter is clear as $X_i,Y_i,X_i^G$ all have the homotopy type of a finite CW-complex without cells of dimension bigger than $N$ by our assumptions. (In order to see why the latter holds for $X_i^G$ recall that $X_i^G$ is a CW-complex which we get by attaching only those $G$-cells of $X_i$ which are of the type $G/G$, as we already remarked before Theorem 2.4.) The map $\pi$ induces a homomorphism of the cohomological long exact sequence of the pair $(Y,X^G)$ into that of the pair $(X,X^G)$. Therefore we have the following commutative ladder of long exact sequences:
$$\CD\rightarrow H^{k-1}(X^G)@>>>H^k(Y,X^G)@>>>
H^k(Y)@>>>H^k(X^G)\rightarrow
\\@VV\text{id}V@VV\pi^*V@VV\pi^*V@V\text{id}VV\\
\rightarrow H^{k-1}(X^G)@>>>H^k(X,X^G)@>>>
H^k(X)@>>>H^k(X^G)\rightarrow\endCD$$
The vertical map $\pi^*:H^k(Y)\rightarrow H^k(X)$ is surjective for every $k\in\Bbb N$ because $H^0(X)=H^0(Y)=\Bbb F_p$ and $H^k(X)=0$ if $k>0$. Hence claim $(ii)$ follows from the five lemma. In particular for every $k\geq0$ the group $G$ acts trivially on $H^k(X,X^G)$ hence
$$E_2^{j,k}\cong H^j(G,\Bbb F_p)
\otimes_{\Bbb F_p}H^k(X,X^G)\cong H^k(X,X^G)\quad(\forall k\in\Bbb N)$$
where we used that $H^j(G,\Bbb F_p)\cong\Bbb F_p$ for every $k\in\Bbb N$ when $\Bbb F_p$ is equipped with the trivial $G$-action. The first part of claim $(iii)$ is now clear. In order to show the second part it will be sufficient to show that for every $i\in I$ every class $c\in H^j(G,H^k(X_i,X^G_i))$ which lies in the image of the map $H^j(G,H^k(Y_i,X^G_i))\rightarrow H^j(G,H^k(X_i,X^G_i))$ induced by $\pi_i^*$ is a permanent cocycle of the spectral sequence
$H^j(G,H^k(X_i,X_i^G))\Rightarrow H^{j+k}(Y_i,X_i^G)$. The latter follows at once from 1.4 of [3] on page 30.\ $\square$
\enddefinition
Let us return to the proof of Theorem 4.5. First note that it is enough to show that $H^k(X,X^G)=0$ for every $k\geq 0$. Indeed in this case we get that $H^k(X^G)=H^k(X)$ for every $k\in\Bbb N$ by looking at the cohomological long exact sequence of the pair $(X,X^G)$. Assume that the claim above is false and let $k\in\Bbb N$ be the largest integer such that $H^k(X,X^G)\neq0$; we know that there is such a number by claim $(i)$ above. By claim $(iii)$ of Lemma 4.6 the terms $E^{j,k}_2$ are formed of permanent cocycles. These cocycles are not coboundaries because the differentials $d_r^{s,t}$ ($r\geq2$) strictly decrease the fibre degree and $E_2$ (and hence $E_r$ for $r\geq2$) does not contain non-zero terms of fibre degree larger than $k$. Hence we have $E^{j,k}_{\infty}=E^{j,k}_2\neq0$ so $H^{j+k}(Y,X^G)\neq0$ for every $j\in\Bbb N$ which contradicts claim $(i)$ of Lemma 4.6.\ $\square$
\definition{Proof of Theorem 2.4} Let $\widetilde X_{\Sigma}$ be the pro-object of $\Cal T_0$ introduced in Definition 3.2. By Proposition 3.3 the latter is cohomologically acyclic mod $p$. Let $s:G\rightarrow\pi^{\Sigma}_1(X,G,x)$ be a section of (2.3.1), and let $h_s$
be continuous left-action of $G$ on $\widetilde X_{\Sigma}$ constructed in Definition 3.7. Let $N$ be a natural number such that $X$ has no cells of dimension bigger than $N$. Because $X_K$ is a finite cover of $X$ the latter is also a finite CW-complex without cells of dimension bigger than $N$ for every $K\in I(X,C_{\Sigma},G)$. Moreover the map $\pi_K:X_K\rightarrow X$ is equivariant with respect to $h^K_s$ and $h$ hence $X_K$ is a $G$-CW complex with respect to the action $h^K_s$ for every $K$ as above. Therefore the action $h_s$ and $\widetilde X_{\Sigma}$ satisfy the assumptions of Theorem 4.5. Hence the subspace $\widetilde X_{\Sigma}^G$ of fixed points of $h_s$ is cohomologically acyclic mod $p$, too. Because $H^0(\widetilde X_{\Sigma}^G)\neq0$ the compact topological space $X_K^G$ is not empty for every $K\in I(X,\Cal C_{\Sigma},G)$. Hence the set $\underline{\widetilde X}_{\Sigma}^G$ is also not empty by Lemma 3.6. By Lemma 3.10 for every point $x\in X^G$ lying in the 
image of the map $\widetilde X_{\Sigma}^G\rightarrow X^G$ the conjugacy class $\iota_X(x)$ contains the section $s$. Therefore the map $\iota^{\Sigma}_{X,0}$ is surjective. Now we only need to show that $\iota^{\Sigma}_{X,0}$ is also injective. Assume that the latter is false. Then there are two connected components $C$ and $C'$ of $X^G$ which both lie in the image of the map $\widetilde X^G_{\Sigma}\rightarrow X^G$ by Lemma 3.9. Let $V\cong\Bbb F_p^2$ be the direct summand of $H^0(X^G)$ corresponding the disjoint union of the two components $C$ and $C'$. Because for every $K\in I(X,\Cal C_{\Sigma},G)$ the image of the map $\pi_K|_{X_K^G}:X_K^G\rightarrow X^G$ contains both $C$ and $C'$ the restriction of the homomorphism $\pi_K^*:H^0(X^G)\rightarrow H^0(X_K^G)$ onto $V$ is injective. By taking the inductive limit we get that the restriction of the homomorphism $H^0(X^G)\rightarrow H^0(\widetilde X^G_{\Sigma})$ induced by the map $\widetilde X^G_{\Sigma}\rightarrow X^G$ onto $V$ is injective, too. But $H^0(\widetilde X^G_{\Sigma})=\Bbb F_p$ by Theorem 4.5. This is a contradiction.\ $\square$
\enddefinition

\heading 5. Applications to geometry 
\endheading

\definition{Definition 5.1} Let $G$ be a group and let $\widehat G$ be its profinite completion. Following Serre (see page 13 of [10]) we will say that $G$ is good if the homomorphism of cohomology groups $H^n(\widehat G, M)\rightarrow H^n(G, M)$ induced by the natural homomorphism $G\rightarrow\widehat G$ is an isomorphism for every finite $G$-module $M$. Let $\Sigma$ be a set of primes and let $p$ be an element of $\Sigma$. We will say that a group $G$ is $(\Sigma,p)$-good if the homomorphism of cohomology groups $H^n(G^{\Sigma}, M)\rightarrow H^n(G, M)$ induced by the natural homomorphism $G\rightarrow G^{\Sigma}$ is an isomorphism for every finite $p$-torsion $G^{\Sigma}$-module $M$. We will say that $G$ is very good if $G$ is $(\Sigma,p)$-good for every prime number $p$ and for every $\Sigma$ containing $p$.
\enddefinition
\proclaim{Lemma 5.2} A group $G$ is $(\Sigma,p)$-good if and only if $\varinjlim_{N\triangleleft G}H^i(N,\Bbb F_p)=0$ for every positive integer $i$ where $N$ ranges over all normal subgroups such that the quotient $G/N$ is in $\Cal C_{\Sigma}$.
\endproclaim
\definition{Proof} Note that the induction of a $p$-torsion module from a subgroup is also $p$-torsion. Hence the argument of Exercise 1 of [10] on page 13 also shows the following: a group $G$ is $(\Sigma,p)$-good if and only if $\varinjlim_{N\triangleleft G}H^i(N,M)=0$ for every positive integer $i$ and for every finite $p$-torsion $G^{\Sigma}$-module $M$, where $N$ ranges over all normal subgroups such that the quotient $G/N$ is in $\Cal C_{\Sigma}$. Therefore the first property in the claim above implies the second. Now assume that $G$ satisfies the second property and let $M$ be a finite $p$-torsion $G^{\Sigma}$-module. We may assume that $M$ is a trivial $G^{\Sigma}$-module by restricting to a normal subgroup $N$ of $G$ such that $G/N$ is in $\Cal C_{\Sigma}$, if this is necessary. In this case $M$ is a direct sum of a finite number of copies of $\Bbb F_p$.
Since cohomology commutes with direct sums in the second variable the claim is now clear.\ $\square$
\enddefinition
\proclaim{Lemma 5.3} Let $\Sigma$ be the set of all primes. Then a group $G$ is good if and only if it is $(\Sigma,p)$-good for every prime number $p$.
\endproclaim
\definition{Proof} The first condition clearly implies the second. Now assume that $G$ is $(\Sigma,p)$-good for every prime number $p$ and let $M$ be a finite $G$-module.  Since cohomology commutes with direct sums in the second variable, we may assume that $M$ is $p$-primary for a prime $p$ while we prove that the homomorphisms $H^n(\widehat G, M)\rightarrow H^n(G,M)$ are isomorphisms. The claim then follows by induction on $n$ and the exponent of $p$ in the prime factorization of the exponent of $G$ by applying the five lemma to the following ladder of cohomological long exact sequences:
$$\longrightarrow H^k(\widehat G,M[p])\longrightarrow H^k(\widehat G,M)\longrightarrow 
H^k(\widehat G,M/M[p])\longrightarrow H^{k+1}(\widehat G,M[p])\longrightarrow$$
$$\CD\quad\quad\quad\quad@VVV
\quad\quad\quad\quad\quad\quad\quad@VVV
\quad\quad\quad\quad\quad\quad\quad@VVV
\quad\quad\quad\quad\quad\quad\quad@VVV
\quad\quad\quad\quad\quad
\endCD$$
$$\longrightarrow H^k(G,M[p])\longrightarrow
H^k(G,M)\longrightarrow
H^k(G,M/M[p])\longrightarrow H^{k+1}(G,M[p])
\longrightarrow$$
where $M[p]$ is the $G$-module of $p$-torsion elements of $M$.\ $\square$
\enddefinition
\definition{Example 5.4} One important consequence of Lemma 5.3 above is that if $G$ is very good then it is also good. Now let $G$ be a finitely generated free abelian group or isomorphic to the fundamental group of a compact Riemann surface. Then for every prime number $p$ the cohomology group $H^*(G,\Bbb F_p)$ is generated by the finite group $H^1(G,\Bbb F_p)$. Hence there is a normal subgroup $N\triangleleft G$ such that $G/N$ is a finite abelian $p$-torsion group and the restriction homomorphism $H^n(G,\Bbb F_p)\rightarrow H^n(N,\Bbb F_p)$ is the zero map for every positive integer $n$. Since every subgroup of $G$ of finite index is a finitely generated free abelian group or isomorphic to the fundamental group of a compact Riemann surface, respectively, the group $G$ is very good by Lemma 5.2. 
\enddefinition
\proclaim{Lemma 5.5} The group $H$ is good if there is a short exact sequence 
$$\CD1@>>>N@>>>H@>>>G@>>>1\endCD$$
such that $G$, $N$ are good, $N$ is finitely generated, and the cohomology groups $H^q(N,M)$ are finite for every $q\in\Bbb N$ and every finite $H$-module $M$.\endproclaim
\definition{Proof} This is the content of Exercise 2(c) of [10] on page 14.\ $\square$
\enddefinition
\definition{Example 5.6} By the lemma above direct products of fundamental groups of compact Riemann surfaces and finitely generated free abelian groups are good. Because finite groups are good, groups which have a good finitely generated normal subgroup $N$ of finite index such that the cohomology groups $H^q(N,M)$ are finite for every $q\in\Bbb N$ and every finite $N$-module $M$ are also good. Moreover extensions of fundamental groups of compact Riemann surfaces by finitely generated free abelian groups are good, too. Finally note that a group $G$ isomorphic to the fundamental group of a non-orientable compact surface $S$ is not very good if the universal cover of $S$ is contractible, but it is good by the above. In fact the completion $G^{\Sigma}$ is the trivial group if $\Sigma$ is the set $\{p\}$  and $p$ is an odd prime, but $H^1(G,M)\neq0$ when $M$ is the induction of the trivial module $\Bbb F_p$ from the two index subgroup of $G$ corresponding to the canonical orientable two-fold cover of $S$.
\enddefinition
The reason why we are interested in this concept is the following
\proclaim{Proposition 5.7} Let $X$ be a connected finite CW-complex whose universal cover is contractible and whose fundamental group is $(\Sigma,p)$-good. Then $X$ is $(\Sigma,p)$-erasable. 
\endproclaim
\definition{Proof} Fix a base point $x\in X$ and by slight abuse of notation let $X$ also denote the pointed space $(X,x)$. Let $\Cal F$ be a $p$-torsion locally constant sheaf with finite stalks on $X$ whose monodromy group is in $\Cal C_{\Sigma}$, and let $M$ be the corresponding $\pi_1(X)$-module. Because the monodromy group of $\Cal F$ is in $\Cal C_{\Sigma}$ the module $M$ is actually the restriction of a $\pi_1(X)^{\Sigma}$-module, which we will denote by the same symbol by the usual abuse of notation. By our assumptions $X_L$ is a finite CW-complex with the homotopy type of the Eilenberg-MacLane space $K(L,1)$ for every $L\in I(X,\Cal C_{\Sigma})$. In particular there is a natural isomorphism:
$$\iota_L:H^k(X_L,\pi_L^*(\Cal F))\rightarrow H^k(L,M).$$
Note that the diagram:
$$\CD H^k(X,\Cal F)@>\iota_{\pi_1(X)}>>H^k(\pi_1(X),M)\\
@V\pi_L^*VV@V\text{Res}VV\\
H^k(X_L,\pi_L^*(\Cal F))@>\iota_L>>H^k(L,M)\endCD$$
is commutative for every $L$ as above. Since $\varinjlim_{L\triangleleft\pi_1(X)}H^k(L,M)=0$ for every positive integer $k$, where $L$ ranges over $I(X,\Cal C_{\Sigma})$, as we saw in the proof of  Lemma 5.2, the claim now follows.\ $\square$
\enddefinition
\definition{Example 5.8} By the proposition above and Example 5.4 compact Riemann surfaces of positive genus and complex tori are $(\Sigma,p)$-erasable for every $\Sigma$ containing $p$. Now let $\Sigma$ be the set of all primes. Then by the proposition above and Example 5.6 a finite topological cover of a product of compact Riemann surfaces of positive genus and complex tori and smooth families of complex tori over hyperbolic compact Riemann surfaces are $(\Sigma,p)$-erasable. For other examples satisfying the conditions of Proposition 5.7 see [7].
\enddefinition
\definition{Proof of Theorem 1.3} By the Grauert-Remmert theorem there is an isomorphism $\pi^{\Sigma}_1(\overline X(\Bbb C),\overline{\eta})\cong\widehat\pi_1^{\Sigma}(\overline X,\overline{\eta})$. For every $K\in I(\overline X(\Bbb C),\Cal C_{\Sigma})$ let $\pi_K:(\overline X_K,\overline{\eta}_K)\rightarrow (\overline X,\overline{\eta})$ be the finite \'etale covering of pointed complex projective varieties corresponding to $K$ under the isomorphism above. For every pair of $K$, $L\in I(\overline X(\Bbb C),\Cal C_{\Sigma})$ such that $L\leq K$ there is a unique map $\pi_{K,L}:(\overline X_L,\overline{\eta}_L)\rightarrow(\overline X_K,\overline{\eta}_L)$ such that $\pi_L=\pi_K\circ\pi_{K,L}$. 
Let $\widetilde X_{\Sigma}$ denote the pro-object of the category $\Cal V_0$ of pointed algebraic varieties over $\Bbb C$ corresponding to the functor $\Cal F_{\Sigma}:\Cal I(\overline X(\Bbb C),\Cal C_{\Sigma})^{\text{op}}\rightarrow\Cal V_0$ which assigns to every object $K$ the pointed space $(\overline X_K,\overline{\eta})$ and to every morphism $\phi_{K,L}$ the map $\pi_{K,L}$. The composition of $\Cal F_{\Sigma}$ with the functor $\Cal P:\Cal V_0\rightarrow\Cal T_0$ which assigns to every pointed algebraic variety $(V,x)$ over $\Bbb C$ the underlying topological space $(V(\Bbb C),x)$ of $\Bbb C$-valued points is the pro-object $\widetilde{\overline X(\Bbb C)}_{\Sigma}$ associated to $\overline X(\Bbb C)$ in Definition 3.2. By definition $\widehat\pi_1^{\Sigma}(X,\overline{\eta})$ is the group of those automorphisms $\gamma$ of the pro-object $\widetilde X_{\Sigma}$ such that $\pi\circ\gamma=\tau\circ\pi$ where $\pi:\widetilde X_{\Sigma}\rightarrow\overline X$ is the natural map and $\tau$ is either the identity or the complex conjugation. Therefore by GAGA the group $\widehat\pi_1^{\Sigma}(X,\overline{\eta})$ is isomorphic via the functor $\Cal P$ to the group $\pi_1^{\Sigma}(\widetilde{\overline X(\Bbb C)}_{\Sigma},G,\overline{\eta})$, where we equip the topological space $\overline X(\Bbb C)=X(\Bbb C)$ with the action of the two-element group $G$ where the generator of $G$ acts via the complex conjugation. Moreover under this isomorphism the short exact sequence (1.1.2) is the same as the short exact sequence (2.3.1) associated to the topological space $\overline X(\Bbb C)$ equipped with the action described above. Because $X(\Bbb R)$ is the set of fixed points of the complex conjugation acting on $X(\Bbb C)$ it is now clear that the claim follows from Lemma 2.2 and Theorem 2.4.\ $\square$
\enddefinition
\definition{Remark 5.9} Recall that a smooth real variety by definition is an ordered pair $(X,\iota)$ where $X$ is a complex manifold and $\iota$ an anti-holomorphic involution of $X$. (Such objects were referred to as real complex manifolds in [8].) Examples are provided by pairs of the form $(X(\Bbb C),\iota)$ where $X$ is a smooth, geometrically connected projective variety defined over $\Bbb R$, and $\iota$ is the complex conjugation. There are other examples however, for example real tori; these are pairs $(X,\iota)$ where $X$ is a quotient of $\Bbb C^n$ by a co-compact discrete lattice $\Lambda\subset\Bbb C^n$ which is invariant under the coordinate-wise complex conjugation and $\iota$ is the unique self-map of $X$ which lifts to the coordinate-wise complex conjugation of $\Bbb C^n$. These are not algebraic if the lattice $\Lambda$ does not admit a Riemann form. One may apply Theorem 2.4 to these objects directly to get a version of Theorem 1.3. As we saw in the proof above such a result is indeed a generalization of Theorem 1.3.
\enddefinition
\definition{Notation 5.10} Let $X$ be a smooth connected projective algebraic curve of positive genus defined over $\Bbb C$ and assume that $G$ is a  finite group such that every non-trivial subgroup of $G$ acts non-trivially on the algebraic curve $X$. Let $\Sigma$ be as above and assume that $G$ is in $\Cal C_{\Sigma}$. Fix an $\eta\in X(\Bbb C)$ and let $\pi^{\Sigma}_1(X(\Bbb C),\eta)\cong\widehat\pi_1^{\Sigma}(X,\eta)$ be the isomorphism furnished by the Grauert-Remmert theorem. For every $K\in I(X(\Bbb C),\Cal C_{\Sigma})$ let $\pi_K:(X_K,\eta_K)\rightarrow (X,\eta)$ be the finite \'etale covering of pointed complex projective varieties corresponding to $K$ under the isomorphism above. For every pair of $K$, $L\in I(X(\Bbb C),\Cal C_{\Sigma})$ such that $L\leq K$ there is a unique map $\pi_{K,L}:(X_L,\eta_L)\rightarrow(X_K,\eta_L)$ such that $\pi_L=\pi_K\circ\pi_{K,L}$. As in the proof above let $\widetilde X_{\Sigma}$ denote the pro-object of the category $\Cal V_0$ introduced above corresponding to the functor $\Cal F_{\Sigma}:\Cal I(X(\Bbb C),\Cal C_{\Sigma})^{\text{op}}\rightarrow\Cal V_0$ which assigns to every object $K$ the pointed space $(X_K,\eta_K)$ and to every morphism $\phi_{K,L}$ the map $\pi_{K,L}$. Let $\widehat\pi^{\Sigma}_1(X,G,\overline{\eta})$ be the group of those automorphisms $\gamma$ of the pro-object $\widetilde X_{\Sigma}$ such that $\pi\circ\gamma=\tau\circ\pi$ where $\pi:\widetilde X_{\Sigma}\rightarrow\overline X$ is the natural map and $\tau\in G$. It is a profinite group and there is a short exact sequence:
$$\CD1@>>>\widehat\pi^{\Sigma}_1(X,\overline{\eta})@>>>
\widehat\pi^{\Sigma}_1(X,G,\overline{\eta})@>>>G@>>>1.\endCD\tag5.10.1$$
Let $\text{Sec}^{\Sigma}(X,G)$ denote the set of conjugacy classes of sections of (5.10.1), defined the same way as in Definition 2.3. Let $X^G$ denote the variety of fixed points of the action of $G$ on $X$.
Every $\Bbb C$-valued point of $X^G$ furnishes a section $G\rightarrow\widehat\pi_1(X,\overline{\eta})$ of the sequence (5.10.1), well-defined up to conjugation. Let
$$\iota^{\Sigma}_{X,G}:X^G(\Bbb C)\rightarrow
\text{Sec}^{\Sigma}(X,G)\tag5.10.2$$
denote the map which sends every $\Bbb C$-valued point of $X^G$ to the corresponding conjugacy class of sections.
\enddefinition
\proclaim{Theorem 5.11} The map $\iota^{\Sigma}_{X,G}$ is a bijection. 
\endproclaim
\definition{Proof} Similarly to Definition 3.5 for every pro-object $Z$ in the category $\Cal V$ of complex algebraic varieties we define its set of points $\underline Z$ to be the set of morphisms from Spec$(\Bbb C)$ to $X$ in the category pro-$\Cal V$. The assignment $X\mapsto\underline X$ is clearly a functor from pro-$\Cal V$ to the category of sets. Note that the theorem above is equivalent to the following claim: for every section $s:G\rightarrow\widehat\pi^{\Sigma}_1(X,G,\overline{\eta})$ there is a unique point $x\in\underline{\widetilde X}_{\Sigma}$ which is fixed by the action of $G$ on $\underline{\widetilde X}_{\Sigma}$ supplied by $s$.

We are going to prove the claim by induction on the order of $G$, that is, we will assume that we proved the claim for every finite group $H$ whose order is strictly less than the order of $G$. Suppose first that $G$ has a proper normal subgroup $H$. Let $\phi:G\rightarrow \widehat\pi^{\Sigma}_1(X,G,\overline{\eta})$ be a section of the short exact sequence (5.10.1). Then there is a commutative diagram:
$$\CD1@>>>\widehat\pi^{\Sigma}_1(X,\overline{\eta})@>>>
\widehat\pi^{\Sigma}_1(X,H,\overline{\eta})@>>>H@>>>1\\
@.@VVV@VVV@VVV@.\\
1@>>>\widehat\pi^{\Sigma}_1(X,\overline{\eta})@>>>
\widehat\pi^{\Sigma}_1(X,G,\overline{\eta})@>>>G@>>>1
\endCD\tag5.11.1$$
because $H$ belongs to the class $\Cal C_{\Sigma}$. The restriction $\phi|_H$ of $\phi$ to $H$ is a section of the short exact sequence in the top line therefore there is a unique fixed point $x\in\underline{\widetilde X}_{\Sigma}$ of the action of $H$ on $\underline{\widetilde X}_{\Sigma}$ corresponding to $\phi|_H$. Because $H$ is a normal subgroup the point $x$ is also fixed by $G$. Hence $\iota^{\Sigma}_{X,G}$ is a bijection.

Assume now that $G$ is a simple group. Then either $G$ is isomorphic to a group of prime order or it is non-abelian.  When $G$ has prime order the theorem follows at once from Theorem 2.4 using the same argument which we used in the proof of Theorem 1.3, and the fact that $X(\Bbb C)$ is $(\Sigma,p)$-erasable (see Example 5.8). Hence we may assume that $G$ is not commutative. If the set $\text{Sec}^{\Sigma}(X,G)$ is empty there is nothing to prove. If it is non-empty then the same also holds for every proper subgroup of $G$, hence the latter have a fixed point by the induction hypothesis. In particular every proper subgroup of $G$ is cyclic. At this point it is possible to derive a contradiction by purely group-theoretical means: let $H$ be a maximal proper subgroup of $G$. For every $x\in G-H$ the subgroup $H\cap x^{-1}Hx$ is normalized by $x$ and it is also normalized by $H$ because the latter is commutative. Hence $H\cap x^{-1}Hx$ is a normal subgroup in $G$ since $x$ and $H$ generate $G$. Therefore it is the trivial subgroup. So $G$ is a Frobenius group. But Frobenius groups are never simple because of the existence of the Frobenius complement.\ $\square$
\enddefinition

\enddocument